\documentclass[11pt]{amsart}
\usepackage{geometry, graphicx, amssymb, amscd, hyperref}
\usepackage[all]{xypic}
\usepackage[latin1]{inputenc}
\usepackage{epstopdf}
\newcommand{\ra}{\rightarrow}
\newcommand{\la}{\leftarrow}
\newcommand{\lla}{\longleftarrow}
\newcommand{\lra}{\longrightarrow}
\newcommand{\hra}{\hookrightarrow}

\newcommand{\rlas}{\rightleftarrows}
\newcommand{\retrmanlabel}[4]{\ensuremath \xymatrix@1{#1\ar@<0.4ex>[r]^-{#3}&#2 \ar@<0.4ex>[l]^-{#4} }}
\newcommand{\retrlabel}[2]{\ensuremath \xymatrix@1{#1\ar@<0.4ex>[r]^-{r}&#2 \ar@<0.4ex>[l]^-{s} }}
\newcommand{\retr}[2]{\ensuremath \xymatrix@1{#1\ar@<0.4ex>[r]&#2 \ar@<0.4ex>[l] }}

\newcommand{\R}{\mathcal{R}}
\newcommand{\Rfd}{\mathcal{R}^{fd}}
\newcommand{\Rld}[1]{\mathcal{R}^{ld}({\mathbb J}#1)}
\newcommand{\Rex}{\mathcal{R}^{\%}}
\newcommand{\V}{\mathcal{V}}
\newcommand{\C}{\mathcal{C}}
\newcommand{\D}{\mathcal{D}}
\newcommand{\simp}[1][\mathfrak B]{{\mathcal S}(#1)}
\newcommand{\Spaces}{{\mathcal S}paces}
\newcommand{\WCat}{{\mathcal WC}at}
\newcommand{\Set}{{\mathcal S}ets}

\newcommand{\holim}[1][ ]{\underset{#1}{\operatorname{holim\ }}}
\newcommand{\hocolim}[1][]{\underset{#1}{\operatorname{ hocolim\ }}}
\newcommand{\holimin}[1][ ]{{\operatorname{holim\ }}_{#1}}
\newcommand{\hocolimin}[1][]{{\operatorname{ hocolim\ }}_{#1}}

\newcommand{\Aex}{A^{\%}}
\newcommand{\wtA}{\widetilde{A}}
\newcommand{\Abiv}{A^{biv}}
\newcommand{\Apara}[2][B]{\ensuremath{A_{#1}\left(#2\right)}}
\newcommand{\Aexpara}[2][B]{\ensuremath{\Aex_{#1}\left(#2\right)}}

\newcommand{\ehom}{\ensuremath\chi^h}
\newcommand{\etop}{\ensuremath\chi^t}
\newcommand{\wtetop}{{\widetilde\chi}^{t}}
\newcommand{\ebiv}{\ensuremath\chi^{biv}}

\newcommand{\thom}{\ensuremath\tau^{h}_{\rho}}
\newcommand{\ttop}{\ensuremath\tau^{t}_{\rho}}
\newcommand{\tsm}{\ensuremath\tau^{s}_{\rho}}

\newcommand{\topm}[1][M]{TOP(#1)}
\newcommand{\tops}{T}

\newcommand{\fib}[1]{\ensuremath{p_{#1}:E_{#1} \to B_{#1}}}
\newcommand{\qfib}{\ensuremath{E \overset{p}{\to} B }}
\newcommand{\utot}[1][M]{E\tops\times_{\tops} #1}

\newcommand{\oM}{\ensuremath\overline M}
\newcommand{\id}{\ensuremath{\rm id}}
\newcommand{\hl}{\ensuremath [0, \infty)}
\newcommand{\J}{\ensuremath {\mathbb J}}
\newcommand{\goth}[1][B]{\mathfrak{#1}}
\newcommand{\tB}[1][p]{#1\mathfrak{B}}
\newcommand{\Map}{{\rm Map}}
\newcommand{\Mor}{{\rm Mor}}

\swapnumbers

\theoremstyle{plain}
\newtheorem{theorem}{Theorem}[section]
\newtheorem{lemma}[theorem]{Lemma}
\newtheorem{proposition}[theorem]{Proposition}
\newtheorem{corollary}[theorem]{Corollary}

\theoremstyle{definition}
\newtheorem{definition}[theorem]{Definition}

\newtheorem{notation}[theorem]{Notation}
\newtheorem{remark}[theorem]{Remark}
\newtheorem{nn}[theorem]{}
\newtheorem*{ack}{Acknowledgments}

\title[Additivity  for topological torsion]{Additivity for the parametrized topological Euler characteristic and Reidemeister torsion}
\date{05/10/2007}                                       

\author[B. Badzioch]{Bernard Badzioch}
\author[W. Dorabia{\l}a]{Wojciech Dorabia{\l}a}

\address[]{Department of Mathematics, SUNY at Buffalo, Buffalo, NY}
\address[]{Department of Mathematics, Penn State Altoona, Altoona, PA }

\begin{document}

\begin{abstract}
Dwyer, Weiss, and Williams have recently defined the notions of the parametrized topological
Euler characteristic and the parametrized topological Reidemeister torsion
which are invariants of bundles of compact topological manifolds.  We show that these invariants 
satisfy additivity formulas paralleling the additive  
properties of the classical Euler characteristic and the Reidemeister 
torsion of $CW$--complexes. 
\end{abstract}

\maketitle

%%%%%%%%%%%%%%%%%%%%%%%%%%%%%%%%%%%%
%
%             INTRODUCTION
%
%%%%%%%%%%%%%%%%%%%%%%%%%%%%%%%%%%%%

\section{Introduction}
\label{INTRO}

In recent years various attempts have been made to generalize the classical Reidemeister
torsion (which is an invariant of non-simply connected finite CW complexes) to a parametrized
version, i.e. to an invariant of fiber bundles. One such generalization is the notion of 
an analytic torsion introduced by Bismut and Lott \cite{Bis-Lott}. It is defined for bundles of smooth 
manifolds satisfying some additional conditions. Another definition was proposed by Igusa and
Klein \cite{Igusa2},\cite{Igusa} who construct parametrized Reidemeister torsion - also for bundles of smooth manifolds - 
using generalized Morse functions. 

Our main interest in this paper lies in yet another definition of torsion developed by
Dwyer, Weiss and Williams \cite{DWW}. One of the main features of their construction is that 
it is described in the language of the homotopy theory which makes it relatively simple. 
Briefly, it proceeds  as follows. Given a manifold $M$ and a locally constant sheaf of $R$-modules $\rho\colon V\to M$
such that the homology groups $H_{\ast}(M, \rho)$ vanish, the classical Reidemeister torsion
of $M$ can be defined as an element of a group ${\rm Wh}(\rho)$ which is a certain quotient of  $K_{1}(R)$. Let 
 $Q(M_{+})$ denote the infinite loops space associated with the suspension spectrum of $M$. 
The starting point for the construction of Dwyer, Weiss and Williams is the observation that  the group
${\rm Wh}(\rho)$ can be identified with $\pi_{0} {\Phi^{s}_{\rho}(M)}$ where 
where $\Phi^{s}_{\rho}(M)$ is  the homotopy fiber of a certain map $\lambda \colon Q(M_{+})\to K(R)$ defined in terms of the sheaf $\rho$.  It follows that the classical torsion of $M$ determines a connectedness component of $\Phi^{s}_{\rho}(M)$. In fact, more is true: one can construct the torsion 
invariant as a unique point $\tau^{s}_{\rho}(M)\in\Phi^{s}_{\rho}(M)$. Assume now that we are given a smooth bundle $p\colon E\to B$ and a sheaf of 
$R$-modules $\rho\colon V\to E$ such that the homology groups of the fibers $M_{b}$ of $p$
with coefficients in $\rho|_{M_{b}}$ vanish.  For every $b\in B$ we obtain a point
$\tau^{s}_{\rho}(M_{b})\in \Phi^{s}_{\rho}(M_{b})$. The spaces $\Phi^{s}_{\rho}(M_{b})$ can be 
glued together using the topology of the bundle $p$ in such way that we obtain a fibration 
$\Phi^{s}_{\rho}(p)\to B$. The assignment $b\mapsto \tau^{s}_{\rho}(M_{b})$ defines then a 
section $\tsm(p)$ of this fibration which can be interpreted as a parametrized smooth torsion of the bundle $p$.

An advantage of the construction of torsion sketched above is its flexibility:  
while $\tsm(p)$ is defined for smooth bundles of manifolds,  its variants can be used 
to extend the notion 
of parametrized torsion to  more general settings. In fact, Dwyer, Weiss and 
Williams gave two additional versions of their definition. The homotopy Reidemeister torsion $\thom(p)$ is defined for any 
fibration $\fib{}$ with homotopy finitely dominated fibers and is obtained just as the smooth torsion, 
but with $Q(M_{+})$ replaced by $A(M)$ -- the Waldhausen $A$-theory of $M$. 
The topological Reidemeister torsion 
$\ttop(p)$ exists whenever $\fib{}$ happens to be a bundle of compact topological manifolds
and is constructed using $\Aex(M)$ - the excisive version of the $A$-theory.

Each of these invariants takes values in a different infinite loop space; thus topological torsion of 
a smooth bundle is not the same as its smooth torsion, but it can be seen as an approximation
of the smooth torsion. Similarly, homotopy torsion can be considered as an approximation 
of the topological torsion when both are defined. 

The work of Goette \cite{GOETTE1},\cite{GOETTE2} and Igusa extended our understanding of the relationship between the analytical torsion 
of Bismut--Lott and the torsion of Igusa--Klein. It is far less clear, however, how these
 notions 
relate to the smooth torsion of Dwyer--Weiss--Williams. Our goal here is to bring these constructions
closer together. The starting point is an axiomatization of parametrized Reidemeister torsion proposed 
by Igusa \cite{AXIOM}. He showed that if a cohomological version of torsion of smooth bundles satisfies 
two conditions, then it is unique up to a scalar multiple. The first condition is the product formula
relating the torsion of a composition of two fibrations $p_{2}\circ p_{1}$ to the torsion of $p_{1}$ and 
$p_{2}$. The second condition is additivity, which describes  the torsion of a pushout of bundles over a space $B$. In \cite{Do} the second author 
showed that a homotopy theoretical analog of additivity is satisfied by the homotopy torsion of 
Dwyer, Weiss and Williams. The main result of our present paper shows that additivity holds 
also for the topological torsion $\ttop$.

\begin{definition}
\label{SPLITBUNDLE}
Let $\fib{}{}$ be a bundle of closed topological manifolds. We say that $p$ admits a fiberwise codimension one splitting if there are subbundles of manifolds 
$p_{i}\colon E_{i}\ra B$ ($i=0,1,2$) such that $E=E_{1}\cup_{E_{0}}E_{2}$, fibers of $p_{1}$, $p_{2}$
are compact submanifolds (with boundary) of the fibers of $p$, and the fibers of $p_{0}$ are the common
boundary of the fibers of $p_{1}$ and $p_{2}$. 
\end{definition}

\begin{theorem}
\label{MAIN}
Suppose that  $\fib{}{}$ is a bundle of closed manifolds which admits a fiberwise 
codimension one splitting into subbundles $p_{i}\colon E_{i}\to B$ for $i=0, 1, 2$. 
Let $R$ be a ring and  $\rho\colon V\to E$ be a locally constant sheaf of 
finitely generated projective left $R$-modules.  Finally, assume that 
$H_{\ast }(p^{-1}(b);\rho)=0$ and $H_{\ast }(p_{i}^{-1}(b);\rho|_{E_i})=0$ for $i=0,1,2$ and all $b \in B$.
Then there exists a preferred homotopy class of paths in the topological Whitehead space 
$\Phi^{t}_{\rho}(p)$ joining $\ttop (p)$ with $j_{1\ast} \ttop(p_1) +j_{2\ast} \ttop(p_2)-j_{0\ast} \ttop(p_0)$, 
where $j_{i\ast}$ is the map induced by the inclusion $j_{i}\colon E_{i}\hookrightarrow E$. 
\end{theorem}

The definitions of $\Phi^{t}_{\rho}(p)$ and $\ttop(p)$ are recalled in Section \ref{DWW}.
We note that this property of $\ttop$ parallels the additivity of the classical combinatorial 
Reidemeister torsion relating torsion of a finite CW--complex $X\cup_{Z}Y$ to the torsions
of $X$, $Y$ and $Z$.  

While the proof of  Theorem \ref{MAIN} is more subtle than in the case of homotopy torsion, 
the  essential idea is to reduce the problem to additivity of $\thom$ and then use the arguments of 
\cite{Do}. We expect that, similarly, a proof of the additivity for smooth torsion may be obtained by reduction to topological case and application of Theorem \ref{MAIN}.  

The main component of the proof of Theorem \ref{MAIN} is the additivity theorem 
for the topological Euler characteristic:

\begin{theorem}
\label{MAIN2}
Let $\fib{}{}$ be a fiber bundle of closed topological manifolds admitting a fiberwise 
codimension one splitting as in Theorem \ref{MAIN}. There exists a preferred homotopy 
class of paths in $\Aex(p)$ 
joining $\etop(p)$ with $j_{1\ast} \etop(p_1) +j_{2\ast} \etop(p_2)-j_{0\ast} \etop(p_0)$.
\end{theorem}

\noindent The space $\Aex(p)$ denotes here the parametrized excisive $A$-theory of $p$, and 
$\etop(p)\in\Aex(p)$ is the topological Euler characteristic of $p$. Again, definitions of these 
notions are sketched in Section \ref{DWW}. 

In the present paper Theorem \ref{MAIN2} serves as a step in establishing 
the additivity of topological torsion, but it has also other potential applications. 
One of the main results of \cite{DWW} says that a bundle 
of compact manifolds $p\colon E\to B$ is fiber homotopy equivalent to a smooth bundle 
if and only if $\etop(p)$ can be lifted to the infinite loop space associated with the parametrized 
suspension spectrum of $p$. From this perspective additivity of the topological Euler 
characteristic provides a tool for computing an obstruction for smoothing of the bundle $p$.

\begin{nn} {\bf Organization of the paper.}
In Section \ref{DWW} we recall the constructions of Dwyer, Weiss, and Williams. In their 
setting the topological (resp. homotopy) Reidemeister torsion can be thought of as a lift 
of the parametrized topological (resp. homotopy) Euler characteristic. In order to prove 
additivity for $\ttop$ it is then enough to verify that additivity holds for the topological Euler 
characteristic, and then show that we can lift the resulting path. As the first step we 
demonstrate (\S \ref{NONPARAM}) that excisive Euler characteristic is additive in 
non-parametrized case, that is for bundles over a one-point space. 
In Section \ref{FLAT BUNDLES}  we show how to extend this result to bundles 
of manifolds with a discrete structure group. 
Subsequently in \S \ref{UNIVERSAL} we show that additivity of the Euler characteristic for arbitrary
bundles follows from additivity for certain universal bundles. Then, in Section \ref{FLAT} 
we show that additivity for the universal bundles follows from additivity of the topological
Euler characteristic for bundles with a discrete structure groups. This completes the proof 
of Theorem \ref{MAIN2}.  Finally, in Section  \ref{ADDITIVITY} we prove Theorem \ref{MAIN}.  
\end{nn}

\begin{ack}
The  authors wish to thank Bruce Williams for conversations which contributed to this 
work. Comments of the referee helped to clarify several passages. 
\end{ack}

%%%%%%%%%%%%%%%%%%%%%%%%%%%%%%%%%%%%
%
%                 DWW CONSTRUCTION
%
%%%%%%%%%%%%%%%%%%%%%%%%%%%%%%%%%%%%

\section{Dwyer-Weiss-Williams constructions}
\label{DWW}

The purpose of this section is to provide a quick review of constructions leading to the definition 
of the topological Reidemeister torsion. We refer to \cite{DWW} for a detailed treatment of this 
subject. We also set here the notation which we will use throughout the paper. In general 
we tried to preserve the notation of \cite{DWW}, although some differences occur. 

\begin{nn}{\noindent\bf Non-parametrized Euler characteristics.}
\label{NONPARA E}
 By a Waldhausen category we will mean a category $\C$ together with a choice of 
two subcategories: a subcategory of cofibrations and a subcategory of weak equivalences satisfying 
the axioms of \cite{Wal1} (in the terminology of \cite{Wal1} such a category $\C$ 
is called a category with cofibrations and weak equivalences). 
Applying the $S_{\bullet}$-construction to $\C$ one obtains $\Omega|wS_{\bullet}\C|$ -- 
the  $K$-theory space of $\C$ \cite[1.3]{Wal0}. This is an infinite loop space which we will denote 
by $K(\C)$.
Every object $c\in \C$ represents a point $[c]\in K(\C)$, and a weak equivalence 
$\varphi\colon c\to c'$ determines a path from $[c]$ to $[c']$. 

 An exact functor of Waldhausen categories $F\colon \C\ra\D$ is a functor preserving 
 the distinguished subcategories and all other relevant structures. Any such functor 
 induces a map of the associated infinite loop spaces $F_{\ast}\colon K(\C)\to K(\D)$. The functor 
$\C\times \C \ra \C$ which assigns to a pair of objects $(c, c')$ their coproduct $c\vee c'$ in $\C$
is exact and defines a map $K(\C)\times K(\C)\ra K(\C)$. This equips  $K(\C)$ with an 
$H$-space structure such that $[c]+[c']=[c\vee c']$. If we have defined 
a suspension functor $\Sigma \colon \C\ra \C$ \cite[p. 349]{Wal0} then the map $K(\C)\ra K(\C)$
induced by $\Sigma$ represents a homotopy inverse with respect to the $H$-space structure 
on $K(\C)$. Thus, for $c\in\C$ we can write $-[c]:= [\Sigma c]$. All Waldhausen categories 
considered here come equipped with suspension functors.

If $c\to c' \to c''$ is a cofibration sequence in $\C$ then there is a path in $K(\C)$ joining $[c']$ 
with $[c]+[c'']$. One way to get such a path is to use Waldhausen's additivity theorem 
\cite[Prop. 1.3.2]{Wal0}. Another way is  to observe that (in the notation of \cite[1.3]{Wal0}) 
a cofibration sequence $c\ra c'\ra c''$ defines a point $[c\ra c'\ra c'']\in|wS_{2}C|$. The restriction 
of the map $|wS_{2}C|\times \Delta^{2}\ra |wS_{\bullet}\C|$ to $[c\ra c'\ra c'']\times \Delta^{2}$ 
yields the desired path \cite[1.3.3]{Wal0}. This second construction is more explicit 
and easily adapts to the parametrized setting (see Lemma \ref{CHAINCOND2}). 
This is the construction we are using throughout the paper. 
\end{nn}

Following \cite[p. 40]{DWW} we will denote by $\Rfd(X)$ the category of homotopy 
finitely dominated retractive spaces over $X$. The objects of $\Rfd(X)$ are 
diagrams 
$$\retrlabel{Y}{X}$$ 
such that  $r\circ s= {\id}_{X}$, $s$ is a cofibration, and $Y$ is a homotopy finitely 
dominated space over $X$. 
The category $\Rfd(X)$ can be equipped with a Waldhausen category structure where 
a morphism in $\Rfd(X)$ is a weak equivalence or a cofibration if its underlying map of spaces
is a homotopy equivalence or, respectively, a map with the homotopy extension property. 
Its  $K$-theory space is denoted $A(X)$ and called the $A$-theory of $X$.

\begin{definition}
\label{ehom}
Let $X$ be a finitely dominated space. The characteristic object $X^{h}\in\Rfd(X)$ is the 
retractive space  
$$\retrlabel{X\times\{-1, 1\}}{X}$$ 
where $s(X)=X\times\{-1\}$, and $r$ is the projection map. 
The homotopy Euler characteristic of $X$ is the point $\ehom(X)\in A(X)$
represented by $X^{h}$.  
\end{definition} 

The assignment $X\mapsto A(X)$ defines a functor on the category of finitely dominated spaces. 
It is not a homology theory since it does not satisfy the excision axiom. By \cite{WW3}
there exists a functor $X\to\Aex(X)$ which in a certain sense is the best possible approximation of 
$A(-)$ by an excisive functor. 
In \cite[\S7]{DWW}  the authors show that if $X$ is an Euclidean neighborhood retract (ENR) 
then $\Aex(X)$ can be  explicitly  constructed using Waldhausen categories. 
We outline this construction next.

For an ENR space $X$ we have a category  $\Rld{X}$ the objects of which 
are diagrams 
$$\retrlabel{Y}{X\times\hl}$$
such that $r\circ s ={\id}_{X\times\hl}$ 
and where $Y$ is homotopy locally finitely dominated as a space over $X\times\hl$ 
\cite[p.48]{DWW}.
Morphisms in $\Rld{X}$ are retractive maps. Given objects
$r_{i}\colon{Y_{i}}\rlas{X\times\hl}\colon s_{i}$ for $i=1, 2$ and morphisms 
$f, g\colon Y_{1}\ra Y_{2}$ we have the notion of controlled homotopy between $f$ and $g$. 
By this we mean a map $H\colon Y_{1} \times [0,1] \ra Y_{2}$ which gives a homotopy between 
$f$ and $g$ in the usual sense and which commutes with the maps $s_{1}, s_{2}$, but commutes 
with the retractions $r_{1}, r_{2}$ only in a relaxed, controlled way \cite[p. 47]{DWW}.   
The category  $\Rld{X}$ can be equipped with a Waldhausen category structure where 
weak equivalences are controlled homotopy equivalences  and  cofibrations are the maps with 
the controlled homotopy extension property.  
The $K$-theory space of $\Rld{X}$ will be denoted by $A^{\J}(X)$. 

We have a functor 
$$I\colon \Rfd(X)\ra \Rld{X}$$
which assigns to a retractive space $r\colon{Y}\rlas{X}\colon s$ over $X$ a retractive space over 
$X\times\hl$
$$\retr{\overline Y}{X\times\hl}$$
where $\overline Y$ is a pushout in the diagram
$$\xymatrix{
X=X\times\{0\} \ar[r]\ar[d]_-{s}& X\times\hl\ar[d] \\
Y\ar[r]& \overline Y \\
}$$ 
The functor  $I$ is an embedding of categories, and -- considered as a functor 
of Waldhausen categories --  it is exact  so it induces a map of the $K$-theory spaces
$$I_{\ast}\colon A(X)\ra A^{\J}(X)$$

Next, let $\V(X)$ denote the category of proper retractive ENRs over $X\times\hl$. This 
is a subcategory of $\Rld{X}$ whose objects $r\colon {Y}\rlas{X\times\hl}\colon s$ satisfy 
the conditions that $Y$ is an ENR and that $r$ is a proper map \cite[7.8]{DWW}. 
Let $J\colon \V(X)\ra \Rld{X}$ be the inclusion functor.  The category $\V(X)$ admits 
a Waldhausen category structure such that $J$ becomes an exact functor \cite[p. 51]{DWW}.
As a consequence we obtain an infinite loop space $V(X):=K(\V(X))$ and a map
$J_{\ast}\colon V(X)\ra A^{\J}(X)$.

\begin{definition}
\label{Aex}
Let $X$ be a compact ENR. The space $\Aex(X)$ is the  homotopy limit
$$\begin{CD}
\Aex(X):=\holim(A(X) @>I_{\ast}>> A^{\J}(X) @<J_{\ast}<<V(X))
\end{CD}$$
We call $\Aex(X)$ the excisive $A$-theory of $X$. The  natural map 
$\alpha\colon\Aex(X)\to A(X)$ is called the assembly map.
\end{definition}

\begin{remark}
\label{Vcontract}
\noindent 1) In \cite[\S7]{DWW} the space $\Aex(X)$ is defined in somewhat different manner, 
as a homotopy fiber of a map $V(X)\ra K(\mathcal{RG}^{ld}(X))$ where $\mathcal{RG}^{ld}(X)$ is the 
Waldhausen category with the same objects as $\Rld{X}$ but with germs 
of retractive maps as morphisms. The resulting infinite loop space is however homotopy equivalent 
to the one described above (see \cite[Lemma 8.7]{DWW}). 

\medskip 

\noindent 2) The space $V(X)$ is in fact contractible \cite[p.52]{DWW}, so $\Aex(X)$ is equivalent to 
a homotopy fiber of the map $I_{\ast}$. The above description of $\Aex(X)$ allows us however 
to perform certain constructions in $\Aex(X)$ combinatorially, on the level of Waldhausen 
categories as follows. For a compact ENR $X$ let $\Rex(X)$ denote the pullback of 
the diagram of categories 
$$\begin{CD}
\Rfd(X) @>I >> \Rld{X} @<J <<\V(X)
\end{CD}$$
Thus, objects of $\Rex(X)$ are pairs $(a, b)$ where $a\in \Rfd(X)$, $b\in \V(X)$, and $I(a)=J(b)$, 
and morphisms are defined similarly.  The category $\Rex(X)$ has the obvious structure 
of a Waldhausen category such that the functors $\Rex(X)\to \Rfd(X)$, $\Rex(X)\to \V(X)$
are exact. It follows that we have a commutative diagram of infinite loop spaces:
$$\xymatrix{
K(\Rex(X)) \ar[r]\ar[d]& V(X)\ar[d] \\
A(X)\ar[r]& A^{\J}(X) \\
}$$ 
As a consequence we obtain a map $K(\Rex(X))\to \Aex(X)$. In particular any object 
$(a, b)\in \Rex(X)$ determines a point $[a,b]\in\Aex(X)$, a weak equivalence 
$(a, b)\to (a', b')$ defines a path from $[a, b]$ to $[a', b']$, and a cofibration sequence
$$(a, b)\to (a', b')\to (a'', b'')$$
in $\Rex(X)$ determines a path in $\Aex(X)$ joining $[a', b']$ with $[a,b]+[a'', b'']$. 
It will be sometimes convenient to describe points and paths in $\Aex(X)$ in this way. 

\end{remark}

We are now ready to define the topological Euler characteristic of a space.

\begin{definition}
\label{etop}
Let $X$ be an ENR. The characteristic object of $X$ in $\V(X)$ is the object $X^{v}$ 
given by the retractive space  
$$\retr{X\times\{0\}\sqcup X\times\hl}{X\times\hl}$$
The topological Euler characteristic of $X$ is the point $\etop(X)\in \Aex(X)$ determined
by the object $X^{t}:=(X^{h}, X^{v})\in\Rex(X)$.  
\end{definition}

Notice that we have $\alpha(\etop(X))=\ehom(X)$.

\begin{remark}
\label{induced maps}
Let $f\colon X\to Y$ be a map of spaces. By abuse of notation by $f_{\ast}$
 we will denote each of the functors induced by $f$: $\Rfd(X)\to \Rfd(Y)$, $\V(X)\to \V(Y)$, 
$\Rld{X}\to \Rld{Y}$, and $\Rex(X)\to \Rex(Y)$, as well as the maps of infinite loop spaces:
$A(X)\to A(Y)$, $V(X)\to V(Y)$, $A^{\J}(X)\to A^{\J}(Y)$, and $\Aex(X)\to \Aex(Y)$. 
\end{remark}

Euler characteristics are not preserved in general by the maps $f_{\ast}$. However, we have 
the following

\begin{lemma}
\label{lax A}
Let $f\colon X\to Y$  be a homotopy equivalence of finitely dominated 
spaces. We have a canonical weak equivalence $f^{h}\colon f_{\ast}(X^{h})\to Y^{h}$
in $\Rfd(Y)$. 
Moreover, if $g\colon Y\to Z$ is another homotopy equivalence, then 
$(gf)^{h}=g^{h}\circ g_{\ast}(f^{h})$.   
\end{lemma}

\begin{proof}
Define 
$$f^{h}:= f\sqcup {\rm id}\colon f_{\ast}(X^{h})=X\sqcup Y\ra Y\sqcup Y$$
The properties of $f^{h}$ are straightforward to check. 
\end{proof}

\begin{corollary}
\label{lax ehom}
If $f\colon X\ra Y$ is a homotopy equivalence of finitely dominated spaces
then there is a canonical path $\omega_{f}$ from $f_{\ast}\ehom(X)$ to $\ehom(Y)$. 
\end{corollary}

\begin{lemma}
\label{lax V}
Let $f\colon X\ra Y$ be a cell-like map \cite{La1}, \cite{La2}, \cite{La3} of compact ENRs. We have a canonical  
weak equivalence 
$f^{v}\colon f_{\ast}(X^{v})\to Y^{v}$ in $\V(X)$. 
Moreover the weak equivalences 
$f^{h}$ and $f^{v}$ satisfy the equations $I(f^{h})=J(f^{v})$, and so they define 
a weak equivalence $f^{t}:=(f^{h}, f^{v}):f_{\ast}(X^{t})\to Y^{t}$ in $\Rex(Y)$.
The weak equivalences  $f^{v}$ and $f^{t}$ satisfy 
a cochain condition analogous to the one described in Lemma \ref{lax A}. 
\end{lemma}

\begin{proof}
The map $f^{v}$ is given by 
$$f^{v}:=f \sqcup {\rm id}
\colon f_{\ast}(X^{v})= X \sqcup Y\times\hl \to Y \sqcup Y\times\hl$$
It is a weak equivalence in $\V(X)$ by \cite[p. 53]{DWW}.
Verification of the remaining properties of $f^{v}$ is straightforward.
\end{proof}

\begin{corollary}
\label{lax etop}
If $f\colon X\ra Y$ is a cell-like map of compact ENRs then there is a canonical path 
$\sigma_{f}$ joining 
$f_{\ast}(\etop(X))$ with $\etop(Y)$. Moreover, if $a\colon \Aex(X)\to A(X)$ is the assembly 
map, then $a(\sigma_{f})=\omega_{f}$ where $\omega_{f}$ is the path from Lemma \ref{lax ehom}.  
\end{corollary}

\noindent We will refer to the properties of $\ehom(-)$ and $\etop(-)$ described in Corollaries
\ref{lax ehom} and \ref{lax etop} as the lax naturality of Euler characteristics.

\begin{nn}
\label{PARAMETRIZATION}
{\noindent \bf Parametrization.} The constructions sketched above can be 
generalized to the setting where the space $X$ is replaced by a 
fibration $\qfib$ with a fiber $F$. Intuitively, one can construct in this case a
fibration $\Apara{E}\to B$  the fiber of which is $A(F)$. An analog of the 
homotopy Euler characteristic in this context is a section of this fibration 
which restricts to $\ehom(F)$ over every point of $B$. A similar idea underlies 
the notions of the parametrized excisive $A$-theory and the parametrized  topological 
Euler characteristic. For technical reasons it is more convenient, however, to 
define the parametrized Euler characteristics in different terms. 
We give these formal definitions first, and then explain how they relate to the above idea. 
\end{nn}

Let $\C$ be a small category, and let $F\colon \C\to\Spaces$ be a functor. Recall 
that $\holimin[\C]F$ is the space of natural transformations
$$\holim[\C]F:=\Map_{\C}(|\C/-|, F)$$
where $|\C/-|$ is the functor which assigns to $c\in\C$ the nerve of the overcategory 
$C/c$. The following fact is implicitly present in \cite{DWW}:

\begin{lemma}
\label{CHAINCOND}
Let $\WCat$ denote the category of Waldhausen categories with exact functors as 
morphisms. Assume that for a functor $F\colon \C\to \WCat$ we have a rule which assigns to 
$c\in\C$ an object $c^{!}\in F(c)$, and to $f\in\Mor_{\C}(c, d)$ a weak equivalence
$f^{!}\colon F(f)(c^{!})\ra d^{!}$ in $F(d)$ in such way that $(g\circ f)^{!}=g^{!}\circ F(f)(f^{!})$.  
Then the assignment $|C/c|\mapsto [c^{!}]$ defines a point  $[c^{!}, f^{!}]\in \holimin[\C]K(F(c))$.
\end{lemma}

\noindent Indeed, if $wF(c)$ denotes the subcategory of weak equivalences in the Waldhausen category 
$F(c)$, then the assignments $c\mapsto c^{!}$, $f\mapsto f^{!}$ as in Lemma \ref{CHAINCOND}
define a point in $\holimin[\C] |wF(c)|$ (where $|wF(c)|$ is the nerve of $wF(c)$). 
Using the natural transformation of functors $|wF(c)|\to K(F(c))$ we obtain a point in 
$\holimin[\C] K(F(c))$ as claimed.

Now, let $\goth\colon\Delta^{op}\to \Set$ be a simplicial set. Denote by $\simp$ 
the category whose objects are all simplices $x\in\goth$, and where morphisms $x\to y$
in $\simp$ come from morphisms $\varphi$ in $\Delta^{op}$ satisfying $\goth(\varphi)(y) =x$ 
(thus $\simp$ is the opposite category of the Grothendieck construction on 
the functor $\goth$ \cite[p.22]{dwyer1}). Let $B$ be the geometric realization of $\goth$, and for 
$x\in\goth$ let $\varphi_{x}\colon\Delta^{|x|}\to B$ denote the characteristic map of $x$. 
Assume that we have a fibration $\qfib$ with a homotopy finitely dominated fiber $F$. We can define a functor 
$$\simp \to \Spaces, \hskip 1cm x\mapsto E_{x}$$ 
where $E_{x}:=\lim(\qfib\overset{\varphi_{x}}{\la}\Delta^{|x|})$. As a result we get a 
functor 
$$F\colon\simp\ra \WCat, \hskip 1cm x\mapsto \Rfd(E_{x})$$

For $x\in\simp$ consider the assignment $x\mapsto E_{x}^{h}$ where $E_{x}^{h}\in \Rfd(X)$
is the characteristic object of $E_{x}$ (\ref{ehom}). Notice that for any morphism
$f\colon x\to y$ in $\simp$ the map $F(f)\colon E_{x}\ra E_{y}$ is a homotopy equivalence, 
so by Lemma \ref{lax A} it defines a weak equivalence 
$F(f)^{h}\colon F(f)_{\ast}(E_{x}^{h}) \ra E_{y}^{h}$ in $\Rfd(E_{y})$ . Lemma \ref{lax A} also
shows that the assignments $x\mapsto E_{x}^{h}$, $f\mapsto F(f)^{h}$ satisfy the conditions
Lemma \ref{CHAINCOND}, and so they  define a point $[E_{x}^{h}, F(f)^{h}]\in\holim[\simp]A(E_{x})$.

\begin{definition}
\label{ehom para}
The homotopy Euler characteristic of a fibration $\qfib$ is the point 
$[E_{x}^{h}, F(f)^{h}]\in\holim[\simp]A(E_{x})$. We will write $\ehom(p):=[E_{x}^{h}, F(f)^{h}]$. 
\end{definition}

Next, notice that the constant maps $A(E_{x})\to \ast$ define a map of homotopy colimits
$$p_{\ast}\colon\hocolim[\simp]A(E_{x})\ra \hocolim[\simp]\ast = B$$
which by \cite[p.180]{Dror} is a quasi-fibration with the fiber $A(F)$. Let $p_{\ast}\colon \Apara{E}\to B$
be the fibration associated to this quasi-fibration.
For all morphisms $x\ra y$ in $\simp$ the map $E_{x}\ra E_{y}$ is a 
homotopy equivalence, and thus so is the map $A(E_{x})\ra A(E_{y})$. 
Therefore as an application of \cite[Prop. 3.12]{dwyer} we obtain 
\begin{proposition} 
\label{sections}
The map
$$\holim[x\in\simp] A(E_{x})=\Map_{\simp}(|\simp/-|, A(E_{(-)}))\ra \Map_{B}(\hocolim[x\in\simp] |\simp/x|, \hocolim[x\in\simp] A(E_{x}))$$
$$ f\mapsto \hocolim f$$
is a weak equivalence. 
\end{proposition}

\noindent Here $\Map_{B}(-, -)$ denotes the mapping space of spaces over $B$. 
Also, in the category of spaces over $B$ we have weak equivalences 
$$\hocolim[x\in\simp] |\simp/x|\simeq B \ \ \ \  {\rm and} \ \ \ \  \hocolim[x\in\simp] A(E_{x})\simeq \Apara{E}$$
which combined with Proposition \ref{sections} give
$$\holim[x\in\simp] A(E_{x})\simeq \Map_{B}(B, \Apara{E})$$
It follows that  $\ehom(p)\in\holim A(E_{x})$ defines a point (unique up to a contractible space 
of choices) in  $\Gamma(p_{\ast})$ -- the space of sections of the fibration 
$p_{\ast}\colon\Apara{E}\to B$. This brings us back to the intuitive construction of $\ehom(p)$
sketched at the beginning of this section.

As we have mentioned above the idea behind the definition of the parametrized topological 
characteristic $\etop(p)$ is similar. The details, however, are more involved. The problem is that 
one cannot (mimicking Definition \ref{ehom para}) define $\etop(p)$ using Lemmas \ref{lax V}
and \ref{CHAINCOND} since the maps $E_{x}\to E_{y}$ are usually not cell-like. In \cite{DWW}
the authors overcome this difficulty by replacing $\goth$ with a new simplicial set $\tB$, and
the functor $F\colon \simp\to \Spaces$ by a new functor $tF\colon \simp[\tB]\ra \Spaces$, 
which sends every morphism in $\simp[\tB]$ to a homeomorphism. Then the assumptions of
Lemma \ref{lax V} are satisfied and we can apply  \ref{CHAINCOND} in order to define a point in 
$\etop(p)\in\holimin[{\simp[\tB]}]  \Aex(tF)$. The details follow. 

\begin{definition}[Compare {\cite[p.13]{DWW}} ]
\label{tB}
Let $\qfib$ be a (locally trivial) fiber bundle, where $B$ is the geometric realization of a simplicial 
set $\goth$.
By $\tB$ we will denote the simplicial set whose $k$-simplices are pairs $(x, \theta)$ where 
$x$ a $k$-simplex in $\goth$, and $\theta$ is an equivalence relation  on $E_{x}$ such 
that the quotient map $E_{x}\to E_{x}^{\theta}$ and the projection $E_{x}\to \Delta^{k}$
give a homeomorphism $E_{x}\to E_{x}^{\theta}\times \Delta^{k}$. 
\end{definition}

Notice that if $F$ is the fiber of the fiber bundle $\qfib$ then for every $(x, \theta)\in \tB$ 
we have $E_{x}^{\theta}\cong F$ and  the assignment $(x, \theta)\mapsto E_{x}^{\theta}$ 
defines a functor $tF\colon \simp[{\tB}]\to\Spaces$
which sends every morphism in $\simp[\tB]$ to a homeomorphism.

\begin{definition}
\label{etop para}
Let $\qfib$ be a fiber bundle whose fiber $F$ is a compact ENR.
The assignments $(x, \theta)\mapsto (E_{x}^{\theta})^{t}$ (see \ref{etop}) , and 
$f\mapsto F(f)^{t}$ (\ref{lax V}),  where $f$ is a morphism in $\tB$, satisfy the 
conditions of Lemma \ref{CHAINCOND}. Therefore they define a point 
$$\etop(p):=[(E_{x}^{\theta})^{t}, F(f)^{t}]\in \holim[{(x, \theta)\in\tB}] \Aex(E_{x}^{\theta})$$
We call $\etop(p)$ the topological Euler characteristic of the bundle $p$. 
\end{definition}

The following fact lets us compare the homotopy and  topological Euler characteristics
of bundles.

\begin{proposition}
\label{tb2b}
Let $\qfib$ be a locally trivial fiber bundle whose fiber is a compact 
topological manifold $M$ (perhaps with boundary). 
We have a commutative diagram 
$$\xymatrix{
\holim[x]\Aex(E_{x}) \ar[r]^{\simeq}\ar[d]_{\alpha}& 
{\holim[(x, \theta)]\Aex(E_{x})}\ar[r]^{\simeq}\ar[d]_{\alpha}&
\holim[(x, \theta)]\Aex(E_{x}^{\theta})\ar[d]^{\alpha} \\
\holim[x] A(E_{x}) \ar[r]^{\simeq}& 
\holim[(x, \theta)] A(E_{x}) \ar[r]^{\simeq}&
\holim[(x, \theta)] A(E_{x}^{\theta}) \\
}$$ 
The vertical maps are induced by assembly maps. 
All horizontal maps are weak equivalences. The horizontal maps on the left are 
induced by the forgetful functor $\simp[\tB]\to \simp$, while the horizontal maps on the right 
come from the natural transformation $E_{x}\mapsto E_{x}^{\theta}$ of functors over $\simp[{\tB}]$.
\end{proposition}

\noindent The maps on the right are weak equivalences by homotopy invariance of homotopy limits. 
The maps on the left are weak equivalences by \cite[Corollary 2.7]{DWW}.

\begin{nn}
\label{FAH}
The implications of Proposition \ref{tb2b} are twofold. On one hand it lets us think about $\etop(p)$ as 
an element of $\holimin[x\in\simp]\Aex(E_{x})$ defined uniquely up to a contractible space of choices. 
Similarly as in \ref{sections} we then get $\holimin[x\in\simp]\Aex(E_{x})\simeq \Gamma(p^{\%}_{\ast})$
where $\Gamma(p^{\%}_{\ast})$ is the space of sections of the fibration $p^{\%}_{\ast}\colon\Aexpara{E}\to B$ associated to the quasi-fibration $\hocolimin[x\in\simp]\Aex(E_{x})\to B$ with the fiber $\Aex(M)$. 
In particular we can think of $\etop(p)$ as a section of $p^{\%}_{\ast}$.
 On the other hand consider the 
image of $\ehom(p)$ under the composition of the bottom maps in the diagram above. 
One can see that the point  $\holimin[{(x, \theta)}]A(E_{x}^{\theta})$ defined in this way can be 
explicitly described by means of Lemma \ref{CHAINCOND} as coming from the assignment 
$(x, \theta)\mapsto (E_{x}^{\theta})^{ah}:=E_{x}\sqcup E_{x}^{\theta}$ for any object 
$(x, \theta)\in\simp[{\tB}]$, and 
$f\mapsto f^{ah}$ for any morphism $f:(x, \theta)\to (y, \theta')$ in $\simp[{\tB}]$
where $f^{ah}$ is the map
$$f^{ah} \colon tF(f)(E_{x}\sqcup E_{x}^{\theta})=E_{y}\sqcup E_{x}^{\theta}\overset{{\rm id}\sqcup tF(f)}{\lra} E_{y}\sqcup E_{y}^{\theta'}$$ 
We can consider the element $[(E_{x}^{\theta})^{ah}, f^{ah}]\in \holimin[(x, \theta)]A(E_{x}^{\theta})$ 
as a re-definition of $\ehom(p)$. Unlike the non-parametrized case this element is not equal of the 
image of $\etop(p)$ under the assembly map. However, the homotopy equivalences 
$E_{x}\sqcup E_{x}^{\theta}\to E_{x}^{\theta}\sqcup E_{x}^{\theta}$ define a canonical path $\sigma_{p}$
in $\holimin[(x, \theta)]A(E_{x}^{\theta})$, joining the image of $\etop(p)$ with the image of $\ehom(p)$.
\end{nn}

\begin{notation}
\label{A(p)}
We will denote by $A(p)$ the homotopy limit $\holimin[{x}]A(E_{x})$. 
Similarly, ${\Aex}(p)$ will denote $\holimin[{(x, \theta)}]\Aex(E_{x}^{\theta})$. 
\end{notation}

\begin{nn}
{\noindent \bf Reidemeister Torsions.}
Let $R$ be a (discrete) ring with identity, and let $Ch^{fd}(R)$ denote the category 
of chain complexes of left projective $R$-modules which are chain homotopy equivalent 
to finitely generated complexes. The category $Ch^{fd}(R)$ can be equipped with a 
Waldhausen category structure where weak equivalences are chain homotopy equivalences 
and cofibrations are chain maps which are split injective on every level. 
The infinite loop space 
$K(Ch^{fd}(R))$ is homotopy equivalent to  the $K$-theory space of the ring $R$ \cite[p. 43]{DWW}.
Thus, from now on we will denote $K(Ch^{fd}(R))$ by $K(R)$.
\end{nn}

Let $\rho\colon V\to E$ be a locally constant sheaf of projective modules.  The sheaf $\rho$ induces a functor $L_{\rho}\colon\Rfd(E)\to Ch^{fd}(R)$ which assigns to any object 
$(X\leftrightarrows E)\in\Rfd(E)$ the  relative singular chain complex $C(X, E, \rho)$ with local coefficients given by  $\rho$. The functor is $L_{\rho}$ is not exact. It is however close enough to 
being exact that it still defines a map  $L_{\rho\ast}\colon A(E)\to K(R)$ if we slightly modify 
the construction of $A(E)$ and $K(R)$ using a variant of the $S_{\bullet}$-construction proposed 
by Thomason (see \cite[p. 43]{DWW}). Assume that $H_{\ast}(E, \rho)=0$. In this case the complex 
$L_{\rho}(E^{h})$ is acyclic, so the map $0\to L_{\rho}(E^{h})$ (where $0$ is the zero chain complex) 
is a weak equivalence in $Ch^{fd}(E)$. This gives a canonical path $\sigma_{\rho}(E)$ in $K(R)$ 
joining $L_{\rho\ast}(\ehom(E))$ with $\ast$--the basepoint of $K(R)$ which is represented by 
the chain complex $0$. 

\begin{definition}
The pair $(\ehom(E), \sigma_{\rho}(E))$ defines a point in $\Phi^{h}_{\rho}(E)$ - the 
homotopy fiber of the map $L_{\rho\ast}$. This point is called the homotopy Reidemeister 
torsion of the space $E$ and is denoted by $\tau^{h}_{\rho}(E)$. 

If $E$ is a compact ENR then let $\Phi^{t}_{\rho}(E)$ denote the homotopy fiber of the map 
$L_{\rho\ast}\alpha$ where $\alpha$ is the assembly map. The point 
$\tau^{t}_{\rho}(E)\in\Phi^{t}_{\rho}(E)$
defined by the pair $(\etop(E), \sigma_{\rho}(E))$ is the topological Reidemeister torsion of $E$.
\end{definition}
\noindent We will call $\Phi^{h}_{\rho}(E)$ and $\Phi^{t}_{\rho}(E)$ the homotopy (resp. topological)
Whitehead spaces.

Next, let  $\qfib$ be a fibration with a homotopy finitely dominated fiber $F$, where as before  $B=|\goth|$, and 
 and let $\rho\colon V\to E$ be a locally constant sheaf of finitely generated projective left $R$-modules. In such case following 
\cite[p. 66]{DWW} we can define fibrations $\Phi_{\rho, B}^{h}(E)\to B$ and 
$\Phi_{\rho, B}^{t}(E)\to B$ with fibers $\Phi^{h}_{\rho}(F)$ and $\Phi^{t}_{\rho}(F)$ respectively. 
If $H_{\ast}(E_{x}, \rho|_{E_{x}})=0$ for all $x\in\goth$ 
these fibrations admit sections which assign $\tau^{h}_{\rho|_{F}}(F)$ (resp. $\tau^{t}_{\rho|_{F}}(F)$) 
to every point $b\in B$.   We can think of these sections as parametrized versions of the Reidemeister torsions. However, similarly as it was the case for the parametrized Euler characteristics  (\ref{PARAMETRIZATION}), this describes torsions only up to a contractible space of choices. 
We will then again need other, more precise definitions.   

For $x\in\simp$ consider the maps $E_{x}\to E$. The induced functors of Waldhausen 
categories $\Rfd(E_{x})\to \Rfd(E)$ define a natural transformation $\eta$ from the functor 
$F\colon \simp\to \WCat$, $F(x)=\Rfd(E_{x})$ to the constant functor over $\simp$ with 
the value $\Rfd(E)$. Therefore we obtain a map
$$\eta_{\ast}\colon A(p)=\holim[x\in\simp]A(E_{x})\to \holim[x\in\simp]A(E)$$
Recall the the homotopy Euler characteristic $\ehom(p)$ was defined as the point 
in $A(p)=\holimin[x] A(E_{x})$ represented by the assignments $x\mapsto E_{x}^{h}$, 
$f\mapsto F(f)^{h}$ (\ref{ehom para}). The image of $\ehom(p)$ under $\eta_{\ast}$ is 
in turn represented by the assignments $x\mapsto E_{x}\sqcup E$, 
$f\mapsto (F(f)\sqcup {\rm id}\colon E_{x}\sqcup E \to E_{y}\sqcup E)$. Assume that 
$H_{\ast}(E_{x}, \rho|_{E_{x}})=0$ for all $x\in\simp$. In this case the relative chain complexes
$C(E_{x}\sqcup E, E, \rho)$ are acyclic, and the weak equivalences $0\to C(E_{x}\sqcup E, E, \rho)$
define a canonical path $\sigma_{\rho}(p)$ joining the basepoint and 
$L_{\rho_{\ast}}\eta_{\ast}(\ehom(p))$ in $\holimin[x] K(R)$. 

\begin{definition}
\label{PARA H TOR}
Let $\Phi^{h}_{\rho}(p)$ be the homotopy fiber of the map 
$L_{\rho\ast}\eta_{\ast}\colon A(p)\to \holimin[x]K(R)$
over the basepoint of $\holimin[x]K(R)$. 
The homotopy Reidemeister torsion of $p$ is the point 
$\tau^{h}_{\rho}(p)\in\Phi^{h}_{\rho}(p)$ given by the pair $(\ehom(p), \sigma_{\rho}(p))$. 
\end{definition}

Assume now that $\qfib$ is a bundle of compact manifolds. Let $\widetilde A^{\%}(p)$
be the pullback of the diagram 
$$\holim[x]A(E_{x})\overset{\simeq}{\lra} \holim[(x, \theta)] A(E_{x}^{\theta})
\overset{\alpha}{\lla}\holim[(x, \theta)]\Aex(E_{x}^{\theta})$$
where the maps are as in Proposition \ref{tb2b}. Recall (\ref{FAH}) that we have a canonical 
path $\sigma_{p}$ in $\holimin[(x, \theta)] A(E_{x}^{\theta})$ joining the images of $\ehom(p)$
and $\etop(p)$.  Thus, the triple $(\etop(p), \sigma_{p}, \ehom(p))$ defines a point  
${\widetilde\chi}^{t}(p)\in\widetilde A^{\%}(p)$.

\begin{definition}
\label{PARA T TOR}
Let $\Phi^{t}_{\rho}(p)$ denote the homotopy fiber of the map 
$$\widetilde A^{\%}(p)\to \holim[x]A(E_{x})\overset{L_{\rho\ast}\eta_{\ast}}{\lra} \holim[x] K(R)$$
If $H_{\ast}(E_{x}, \rho|_{E_{x}})=0$ for all $x\in\simp$ then the pair 
$({\widetilde\chi}^{t}(p), \sigma_{\rho}(p))$ (where $\sigma_{\rho}(p)$ is the path as in 
Definition \ref{PARA H TOR}) defines a point $\ttop(p)\in\Phi^{t}_{\rho}(p)$. We
call it the topological parametrized Reidemeister torsion of the bundle $p$. 
\end{definition}

Notice that we have a pullback diagram 
$$\xymatrix{
\Phi^{t}_{\rho}(p) \ar[r]^{\gamma_{2}}\ar[d]_{\gamma_{1}}& \wtA^{\%}(p)\ar[d]^{\delta_{1}} \\
\Phi^{h}_{\rho}(p)\ar[r]_{\delta_{2}}& A(p) \\
}$$ 
and that $\gamma_{1}(\ttop(p))=\thom(p)$, $\gamma_{2}(\ttop(p))=\wtetop(p)$.

%%%%%%%%%%%%%%%%%%%%%%%%%%%%%%%%%%%%
%
%                 NON-PARAMETRIZED ADDITIVITY
%
%%%%%%%%%%%%%%%%%%%%%%%%%%%%%%%%%%%%

\section{Non-parametrized additivity theorem}
\label{NONPARAM}

The goal of this section is to prove the following

\begin{theorem}[Additivity for the topological Euler characteristic]
\label{nonparam top additivity}
Let $M$ be a closed topological manifold which admits a splitting along a 
compact codimension one submanifold $M_{0}$:
$$M= M_{1}\cup_{M_{0}}M_{2}$$
There exists a preferred path $\omega$ in $\Aex(M)$ from $\etop(M)$ to 
$k_{1\ast}\etop M_{1}+k_{2\ast}\etop M_{2} - k_{0\ast}\etop M_{0}$, 
where $k_{i}\colon M_{i}\hra M$ is the inclusion map ($i=0, 1, 2$).
\end{theorem}

\begin{figure}[h] 
   \centering
   \includegraphics[width=6cm]{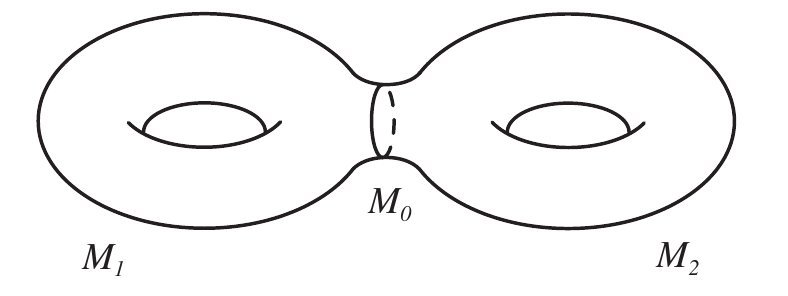} 
\end{figure}

Applying the assembly map $\alpha\colon\Aex(M)\ra A(M)$ to the path $\omega$  we obtain

\begin{corollary} 
\label{nonparam hom additivity}
If $M$ is a manifold as in Theorem \ref{nonparam top additivity}, then there 
exists a path in $A(M)$ joining $\ehom(M)$ with 
$k_{1\ast}\ehom M_{1}+k_{2\ast}\ehom M_{2} - k_{0\ast}\ehom M_{0}$. 
\end{corollary}

\noindent Thus we recover the additivity theorem for the homotopy Euler characteristic 
which was proved in \cite{Do} by the second author. 

\begin{proof}[Proof of Theorem \ref{nonparam top additivity}]
Consider a manifold 
$$\oM=M_{1}\cup_{M_{0}\times \{-1\}}M_{0}\times [-1, 1]\cup_{M_{0}\times\{1\}}M_{2}$$ 

\begin{figure}[h] 
   \centering
   \includegraphics[width=7.5cm]{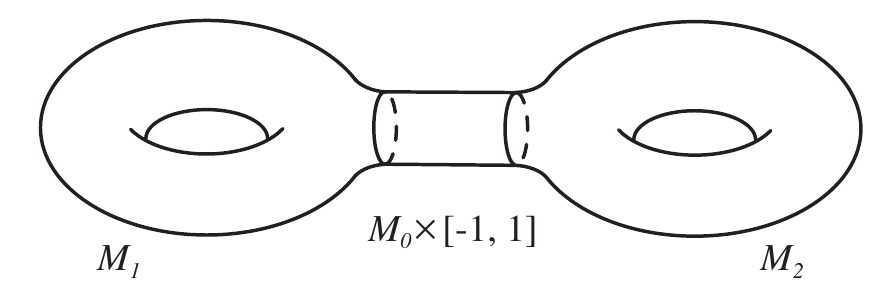} 
\end{figure}

\noindent For  $i=1, 2$ let ${\bar k}_{i}\colon M_{i}\hra\oM$ denote the inclusion map. We have a map 
$f\colon \oM \to M$ which restricts to the inclusions $k_{1}$, $k_{2}$ on $M_{1}$, $M_{2}$
respectively, and which sends $M_{0}\times [-1, 1]\subseteq \oM$ to $M_{0}\subseteq M$
via the projection map onto the first factor. 

We will construct the path $\omega$ as a concatenation of three paths in $\Aex(M)$:
\begin{itemize}
\item[1)] a path $\omega_{1}$ from $\etop(M)$ to $f_{\ast}\etop(\oM)$,
\item[2)] a path $\omega_{2}$ from $f_{\ast}\etop(\oM)$ to 
$f_{\ast}({\bar k}_{1\ast}\etop(M_{1})+{\bar k}_{2\ast}\etop(M_{1})+C)$ for some $C\in \Aex(\oM)$,  
\item[3)] a path $\omega_{3}$ from 
$f_{\ast}{\bar k}_{1\ast}\etop(M_{1})+f_{\ast}{\bar k}_{2\ast}\etop(M_{1})+ f_{\ast}C$ to 
$k_{1\ast}\ehom M_{1}+k_{2\ast}\ehom M_{2} - k_{0\ast}\ehom M_{0}$. 
\end{itemize}

\vskip .2cm
\noindent{\it Construction of $\omega_{1}$.} 
Since $f$ is a cell-like map, the path $\omega_{1}$ exists by lax naturality of $\etop$ (\ref{lax etop}).
 
 \vskip .2cm
\noindent{\it Construction of $\omega_{2}$.}
We will construct a point $C\in \Aex(\oM)$ and a path $\sigma$ in $\Aex(\oM)$ joining 
$\etop(\oM)$ with $({\bar k}_{1\ast}\etop(M_{1})+{\bar k}_{2\ast}\etop(M_{1})-C$. Then we will have 
$\omega_{2}=f_{\ast}(\sigma)$. 

First we define the point $C$ as follows. Let $C' \in \Rfd(\oM)$ be the retractive space over $\oM$
given by 
$$C':=\oM\cup_{M_{0}\times\{-1, 1\}}M_{0}\times[-1, 1]$$ 
and let $C''\in \V(\oM)$ be given by the retractive space over $\oM\times\hl$
$$C'':=\oM\times[0, \infty)\cup_{M_{0}\times\{-1, 1\}}M_{0}\times[-1, 1]$$
where the embedding $M_{0}\times\{-1, 1\}\hra \oM\times[0, \infty)$ is given by
$(x, \pm 1)\mapsto ((x, \pm 1), 0)$. One can check that (in the notation of Remark \ref{Vcontract})
we have $I(C')=J(C'')$, so $(C', C'')$ is an object of $\Rex(\oM)$. Take $C$ to be the point 
in $\Aex(\oM)$ represented by this object.

Next, recall (\ref{etop}) that for any compact ENR $X$ the Euler characteristic $\etop(X)$ is represented 
by the object $X^{t}=(X^{h}, X^{v})\in \Rex(X)$ . Also, notice that in our case the object 
${\bar k}_{i\ast}M_{i}^{h}\in \Rfd(\oM)$ is the retractive space $M_{i}\sqcup \oM$ over $\oM$. 
We have a cofibration sequence in $\Rfd(\oM)$
$${\bar k}_{1\ast}M_{1}^{h}\to (M_{1}\sqcup M_{2}\sqcup \oM) \to {\bar k}_{2\ast}M^{h}_{2}$$
Similarly,  the object ${\bar k}_{i\ast}M^{v}_{i}\in \V(\oM)$ is given by the retractive space 
$M_{i}\sqcup \oM\times\hl$ over $\oM\times\hl$, which gives a cofibration sequence 
in $\V(\oM)$
$${\bar k}_{1\ast}M_{1}^{v}\to (M_{1}\sqcup M_{2}\sqcup \oM\times\hl) \to {\bar k}_{2\ast}M^{v}_{2}$$
These two cofibration sequences lift to a cofibration sequence in $\Rex(M)$:
$${\bar k}_{1\ast}M_{1}^{t}\to (M_{1}\sqcup M_{2}\sqcup \oM\ , \ M_{1}\sqcup M_{2}\sqcup \oM\times\hl) \to {\bar k}_{2\ast}M^{t}_{2}$$
which gives us a path in $\Aex(\oM)$ joining ${\bar k}_{1\ast}\etop(M_{1})+{\bar k}_{2\ast}\etop(M_{2})$
with the point $[M_{1}\sqcup M_{2}\sqcup \oM\ , \ M_{1}\sqcup M_{2}\sqcup \oM\times\hl]$.
As a consequence we only need to construct a path from $\etop(\oM)$ to 
$[M_{1}\sqcup M_{2}\sqcup \oM\ , \ M_{1}\sqcup M_{2}\sqcup \oM\times\hl]+C$. In order to accomplish 
this notice that $C'$ fits into a cofibration sequence in $\Rfd(\oM)$
\begin{equation*}
\label{C'}
\begin{CD}
(M_{1}\sqcup M_{2})\sqcup \oM @>({\bar k}_{1}\sqcup {\bar k}_{2})\sqcup 1>> \oM^{h}=\oM\sqcup\oM
\lra C'
\end{CD}\end{equation*}
while $C''$ is the cofiber in the following cofibration sequence in $\V(\oM)$:
\begin{equation*}
\label{C''}
\begin{CD}
(M_{1}\sqcup M_{2})\sqcup \oM\times[0, \infty) @>({\bar k}_{1}\sqcup {\bar k}_{2})\sqcup 1>> 
\oM^{v}=\oM\sqcup\oM\times[0, \infty)\lra C''
\end{CD}
\end{equation*}
As before these sequences yield a cofibration sequence in $\Rex(\oM)$
$$(M_{1}\sqcup M_{2}\sqcup \oM \ ,\  M_{1}\sqcup M_{2}\sqcup \oM\times\hl)\to 
\oM^{t}=(\oM^{h}, \oM^{v})\to (C', C'')$$
Passing from $\Rex(\oM)$ to $\Aex(\oM)$ we obtain the desired path.

\vskip .2cm
\noindent{\it Construction of $\omega_{3}$.}
We will show that in $\Aex(M)$ we have paths $\delta_{i}$ ($i=1, 2$) joining $f_{\ast}{\bar k}_{i\ast}\etop M_{i}$
with $k_{i\ast}\etop M_{i}$, and a path $\delta_{0}$ from $f_{\ast}C$ (where $C\in \Aex(\oM)$ is defined 
as above) to $-k_{0\ast}\etop M_{0}$. Then we can take $\omega_{3}=\delta_{1}+\delta_{2}+\delta_{0}$. 

For $i=1, 2$ we have $f\circ {\bar k}_{i}=k_{i}$, so $f_{\ast}{\bar k}_{i\ast}\etop M_{i}=k_{i\ast}\etop M_{i}$, and we 
can choose $\delta_{1}, \delta_{2}$ to be the constant paths.

The construction of the path $\delta_{0}$ resembles the construction of $\omega_{2}$ above. 
We have $f_{\ast}C = ([f_{\ast}C'], [f_{\ast}C''])$. Notice that $f_{\ast}C'$ fits into 
the following pushout diagram in $\Rfd(M)$:
$$\xymatrix{
M_{0}\times\{-1, 1\} \ar[r]\ar[d]& M\times[-1, 1]\ar[d] \\
M\ar[r]& f_{\ast}C' \\
}$$ 
Similarly, $f_{\ast}C''$ can be represented as a pushout in $\V(M)$:
$$\xymatrix{
M_{0}\times\{-1, 1\} \ar[r]\ar[d]& M\times[-1, 1]\ar[d] \\
M\times[0, \infty)\ar[r]& f_{\ast}C'' \\
}$$ 
where the map $M_{0}\times\{-1, 1\} \to M\times[0, \infty)$ is given by 
$(x, \pm 1)\mapsto (x, 0)$.
As a consequence we have a cofibration sequence in the category $\Rfd(M)$:
\begin{equation*}
\label{fC'}
\begin{CD}
M_{0}\sqcup M \lra M_{0}\times[-1, 1]\cup_{M_{0}\times\{-1\}}M 
\lra f_{\ast}C'
\end{CD}
\end{equation*}
as well as a cofibration sequence in $\V(M)$:
\begin{equation*}
\label{fC''}
\begin{CD}
M_{0}\sqcup M\times[0, \infty) \lra M_{0}\times[-1, 1]\cup_{M_{0}\times\{-1\}}M\times[0, \infty) 
\lra f_{\ast}C''
\end{CD}
\end{equation*}
(we identify here $M_{0}\times\{-1\}$ with a subspace of $M\times[0, \infty)$ via the embedding
$(x, -1)\mapsto (x, 0)$). 
Notice that  $M_{0}\sqcup M=k_{0\ast}M_{0}^{h}$ in $\Rfd(M)$, and 
$M_{0}\sqcup M\times[0, \infty)= k_{0\ast}M_{0}^{v}$ in $\V(M)$.  
Again, we can lift these two sequence to a cofibration in $\Rex(M)$
$$k_{0\ast}(M^{h}, M^{v})\to 
(M_{0}\times[-1, 1]\cup_{M_{0}\times\{-1\}}M\ ,\  M_{0}\times[-1, 1]\cup_{M_{0}\times\{-1\}}M\times[0, \infty)) 
\to f_{\ast}(C', C'')$$
As a consequence we obtain a path $\delta_{0}'$ in $\Aex(M)$ joining $k_{0\ast}\etop(M_{0}) + f_{\ast}C$
with the point
$$[M_{0}\times[-1, 1]\cup_{M_{0}\times\{-1\}}M\ , \ M_{0}\times[-1, 1]\cup_{M_{0}\times\{-1\}}M\times[0, \infty)]\in \Aex(M)$$
Finally notice that the retractive spaces $M_{0}\times[-1, 1]\cup_{M_{0}\times\{-1\}}M$
and $M_{0}\times[-1, 1]\cup_{M_{0}\times\{-1\}}M\times[0, \infty)$ are weakly equivalent 
to the trivial retractive spaces $M\leftrightarrows M$, and (respectively) 
$M\times[0, \infty)\leftrightarrows M\times[0, \infty)$, which implies that  $\delta_{0}'$
can be further extended to the basepoint of $\Aex(M)$. 
The path $\delta_{0}$ can be now 
obtained shifting $\delta_{0}'$ by the element  $-k_{0}\etop(M_{0})\in \Aex(M)$.

\end{proof}

%%%%%%%%%%%%%%%%%%%%%%%%%%%%%%%%%%%%
%
%                 FLAT BUNDLES
%
%%%%%%%%%%%%%%%%%%%%%%%%%%%%%%%%%%%%
\section{Flat bundles}
\label{FLAT BUNDLES}

The argument which led us to the proof of Theorem \ref{nonparam top additivity}
can be generalized to give a proof of Theorem \ref{MAIN2} for a certain class of fiber 
bundles. Namely, assume that $M$ is a closed topological manifold with a codimension one splitting 
$M= M_{1}\cup_{M_{0}}M_{2}$, and let $G$ be a discrete group acting on $M$ on the right by
homeomorphisms  which preserve the splitting. In such case the bundle 
$p^{G}\colon EG\times_{G}M\to BG$ splits into subbundles $p^{G}_{i}\colon EG\times_{G}M_{i}\to BG$.
For $i=0, 1, 2$ let $j_{i}\colon EG\times_{G}M_{i}\to EG\times_{G}M$ denote the inclusion  map. 
We have 

\begin{proposition}
\label{DISCR BUNDLE ADD}
For the bundle $p^{G}$ as above we have a path in $\Aex(p^{G})$ joining 
$\etop(p^{G})$ with $j_{i\ast}\etop(p^{G}_{1})+j_{i\ast}\etop(p^{G}_{1})-j_{i\ast}\etop(p^{G}_{1})$.
\end{proposition}

The proof of this fact will rely on two lemmas. The first lemma generalizes lax naturality
of the topological Euler characteristic (\ref{lax etop}).  Recall (\ref{lax V}) that 
if $\C$ is a small category and if $F\colon \C\to\Spaces$ is a diagram of compact ENRs then 
we have  assignments $c\mapsto F(c)^{t}$ for $c\in\C$ and $f\mapsto F(f)^{t}$ for a morphism $f$
in $\C$, which  by Lemma \ref{CHAINCOND} define a point $[F(c)^{t}, F(f)^{t}]\in\holimin[c\in\C]\Aex(F(c))$.

\begin{lemma}
\label{LAX HOLIM}
Let $\C$ be a small category. Assume that $F, G\colon \C\to \Spaces$ are diagrams 
of compact ENRs and cell-like maps, and let  $\eta\colon F\to G$ be a natural 
transformation such that for every $c\in \C$ the map $\eta_{c}\colon F(c)\to G(c)$ is cell-like. 
Then there is a path in $\holimin[c\in\C]\Aex(G(c))$ joining $\eta_{\ast}[F(c)^{t}, F(f)^{t}]$
with $[G(c)^{t}, G(f)^{t}]$. Here 
$$\eta_{\ast}\colon\holimin[c\in\C]\Aex(F(c))\to \holimin[c\in\C]\Aex(G(c))$$
is the map induced by $\eta$.
\end{lemma}

Proof of this fact resembles justification for Lemma \ref{CHAINCOND}. The natural 
transformation $\eta$ defines a path in $\holimin[c\in\C]|w\Rex(G(c))|$. Using the map 
$\holimin[c\in\C]|w\Rex(G(c))|\to \holimin[c\in\C]\Aex(G(c))$ we obtain the required path.  

The second lemma describes  how one can construct paths in homotopy limits of diagrams 
of $K$-theory spaces using cofibration sequences.

\begin{lemma}
\label{CHAINCOND2}
Let $F\colon \C\to \WCat$ be a functor, and for $i=1, 2, 3$ let $c\mapsto c_{i}^{!}$, 
$f\mapsto f_{i}^{!}$ be assignments as in the Lemma \ref{CHAINCOND}. Assume also
that for every $c\in\C$ we have a cofibration sequence in $F(c)$:
$$c_{1}^{!}\overset{\varphi_{c}}{\to}c_{2}^{!}\overset{\phi_{c}}{\to}c_{3}^{!}$$
such that for any morphism $c\to d$ in $\C$ the following diagram commutes:
$$ \xymatrix{
F(f)(c_{1}^{!})\ar[d]_{f_{1}^{!}}\ar[r]^{F(f)\varphi_{c}}
& F(f)(c_{2}^{!})\ar[d]^{f_{2}^{!}}\ar[r]^{F(f)\phi_{c}}
&  F(f)(c_{3}^{!})\ar[d]_{f_{3}^{!}}\\
d_{1}^{!}\ar[r]_{\varphi_{d}}
& d_{2}^{!}\ar[r]_{\phi_{d}}
& d_{3}^{!} \\
}
$$ 
Then there is a path in $\holimin[\C]K(F(c))$ joining the point $[c_{2}^{!}, f_{2}^{!}]$
with $[c_{1}^{!}, f_{1}^{!}]+[c^{!}_{3}, f_{3}^{!}]$. 
\end{lemma}

Indeed, in the notation of \cite{Wal1} the cofibration sequences 
$c_{1}^{!}\overset{\varphi_{c}}{\to}c_{2}^{!}\overset{\phi_{c}}{\to}c_{3}^{!}$
define a point 
$[c_{1}^{!}\overset{\varphi_{c}}{\to}c_{2}^{!}\overset{\phi_{c}}{\to}c_{3}^{!}]\in 
\holimin[c\in\C]|wS_{2}F(c)|$. 
Also, we have a map
$$\holim[c\in\C]|wS_{2}F(c)|\times \Delta^{2}\to \holim[c\in\C] |wS_{\bullet}F(c)|$$ 
Restricting this map to 
$[c_{1}^{!}\overset{\varphi_{c}}{\to}c_{2}^{!}\overset{\phi_{c}}{\to}c_{3}^{!}]\times\Delta^{2}$
we obtain the desired path. 

\begin{proof}[Proof of Proposition {\ref{DISCR BUNDLE ADD}}]
Consider the group $G$ as a category with one object $\ast$. The action of $G$ on $M$
defines functors  $F\colon G^{op}\to Spaces$ and $F_{i}\colon G^{op} \to Spaces$ 
(where $G^{op}$ is the opposite category of $G$) such that 
$F(\ast)=\Aex(M)$ and $F_{i}(\ast)=\Aex(M_{i})$
for $i=0, 1, 2$. One can check that we have weak 
equivalences.
$$\Aex(p^{G})\simeq \holim[G^{op}]\Aex(M) \ \ \  {\rm and} \ \ \  \Aex(p_{i}^{G})\simeq \holim[G^{op}]\Aex(M_{i})$$
Moreover, the maps 
$$k_{i\ast}\colon\holim[G^{op}]\Aex(M_{i})\to \holim[G^{op}]\Aex(M)$$ 
induced by the 
inclusions $k_{i}\colon M_{i}\hra M$ correspond under these weak equivalences to the maps 
$j_{i\ast}\colon\Aex(p^{G}_{i})\to \Aex(p^{G})$. Also, the point of $\holimin[G^{op}]\Aex(M)$
corresponding to $\etop(p^{G})$ is (in the notation of  \ref{etop} and \ref{lax V})
the point $[M^{t}, F(g)^{t}]$, and similarly the Euler characteristics $\etop(p^{G}_{i})$, $i=0, 1, 2$
correspond to $[M_{i}^{t}, F_{i}(g)^{t}]\in \holimin[G^{op}]\Aex(M_{i})$. 
As a consequence it is enough to construct a path $\omega$ in $\holimin[G^{op}]\Aex(M)$ joining 
the point $[M^{t}, F(g)^{t}]$ with 
$k_{1\ast}[M_{1}^{t}, F(g)^{t}]+k_{2\ast}[M_{2}^{t}, F(g)^{t}]- k_{0\ast}[M_{0}^{t}, F(g)^{t}]$.

The construction of this path follows the same steps as the proof of 
Theorem \ref{nonparam top additivity}. As in that proof we construct a manifold 
$\oM=M_{1}\cup_{M_{0}\times \{-1\}}M_{0}\times [-1, 1]\cup_{M_{0}\times\{1\}}M_{2}$.
Any homeomorphism of $M$ extends to a homeomorphism of $\oM$ (which is a 
product map on $M_{0}\times [-1, 1]$), thus the group $G$ acts on $\oM$, and we have a 
functor ${\bar F}\colon G \to \Rex(\oM)$ such that ${\bar F}(\ast)=\Rex(\oM)$.  
The map $f\colon \oM\to M$ contracting $M_{0}\times [-1, 1]\subset\oM$ to $M_{0}\subset M$
is $G$-equivariant, so it defines a  natural transformation $f_{\ast}\colon {\bar F}\to F$
which in turn induces a map of homotopy limits
$$f_{\ast}\colon \holim[G^{op}]\Aex(\oM)\to \holim[G^{op}]\Aex(M)$$
Similarly as in the proof of Theorem \ref{nonparam top additivity} we can now build 
the path $\omega$ as a concatenation of three paths. The first path in \ref{nonparam top additivity}, 
(joining $f_{\ast}\etop(\oM)$ with $\etop(M)$) was obtained using lax naturality of $\etop$. 
In our present context we can construct a path in $\holimin[G^{op}]\Aex(M)$ joining 
$f_{\ast}[\oM^{t}, {\bar F}(g)^{t}]$ with $[M^{t}, F(g)^{t}]$ by applying Lemma 
\ref{LAX HOLIM} to the natural transformation $f_{\ast}\colon {\bar F}\to F$. 
The other two paths in the proof of Theorem \ref{nonparam top additivity}
were constructed using certain cofibration sequences in $\Rex(\oM)$ and $\Rex(M)$. 
One can start with the same cofibration sequences and then for each retractive
space appearing in them choose assignments $\{g\to g^{!}\}_{g\in G}$ in such way 
that  the conditions of Lemma \ref{CHAINCOND2} are satisfied. As a consequence 
one obtains the desired paths in $\holimin[G^{op}]\Aex(\oM)$ and $\holimin[G^{op}]\Aex(M)$.
\end{proof}

%%%%%%%%%%%%%%%%%%%%%%%%%%%%%%%%%%%%
%
%                 UNIVERSAL BUNDLES
%
%%%%%%%%%%%%%%%%%%%%%%%%%%%%%%%%%%%%

\section{Universal bundles}
\label{UNIVERSAL}

In the last section we verified that Theorem \ref{MAIN2} holds for some bundles  whose 
structure group is discrete. Our strategy of proving Theorem \ref{MAIN2}
in its whole generality is as follows. In this section we show that the bundles which admit 
the required splitting into subbundles are induced by some universal bundle $p^{U}$. 
Moreover, additivity of the topological Euler characteristic for an arbitrary bundle follows 
from its additivity  for this universal bundle. Thus, we only need to show that $\etop$ 
is additive for the bundle $p^{U}$.  We prove this special case, using the results of the last section, 
in  \S \ref{FLAT}. 

Let $M$ be a closed topological manifold which admits a splitting along 
a codimension one submanifold $M_{0}$:
$$M\simeq M_{1}\cup_{M_{0}}M_{2}$$
Let $\topm$ be the simplicial group of homeomorphisms of $M$ and 
let $\tops$ denote the subgroup of $\topm$ consisting of homeomorphisms
which preserve the splitting. Consider the bundle 

$$p^{U}\colon \utot\ra B\tops$$
The bundle $p^{U}$ admits a fiberwise codimension one splitting into sub-bundles 
$p^{U}_{i}\colon \utot[M_{i}]\ra B\tops$ for $i=0, 1, 2$, such that the fiber 
of $p^{U}_{i}$ is $M_{i}$. Moreover, $p^{U}$ is the universal bundle 
for  bundles $\qfib$ with fiber $M$ which admit such a splitting. Thus, if $p$
splits into subbundles   $p_{i}$ with fibers $M_{i}$ for $i=0, 1, 2$
then we have a map $c\colon B\ra B\tops$ which fits into
pullback diagrams:

\begin{equation}
\label{UPULL}
\xymatrix{
E \ar[r]^-{\bar{c}}\ar[d]_{p}& \utot[M]\ar[d]^{p^{U}} \\
B\ar[r]^{c} & B\tops \\
}
\hskip 2cm 
\xymatrix{
E_{i} \ar[r]^-{\bar{c}_{i}}\ar[d]_{p_{i}}& \utot[M_{i}]\ar[d]^{p^{U}_{i}} \\
B\ar[r]^{c} & B\tops \\
}
\end{equation}
Moreover, these maps commute with all inclusions of sub-bundles. 
Our goal will be to show that the statement of Theorem \ref{MAIN2} holds 
for the bundle $p^{U}$.

\begin{proposition}
\label{universal add}
Let $\iota_{i}\colon  \utot[M_{i}]\hra\utot$ be the inclusion map ($i=0, 1, 2$). 
Then there is a path $\sigma_{p^{U}}$ in $\Aex(p^{U})$ joining $\etop(p^{U})$ with 
$\iota_{1\ast}\etop(p^{U}_{1})+\iota_{2\ast}\etop(p^{U}_{2})-\iota_{0\ast}\etop(p^{U}_{0})$. 
\end{proposition}

\noindent  We postpone the proof of this fact until \S \ref{FLAT}. Meanwhile we will
show that,   as indicated at the beginning of this section,
Theorem  \ref{MAIN2} can be obtained from this special case. We will need the  following 
lemmas which follow directly from the constructions of $\Aex(p)$,  $\etop(p)$, and from 
Lemma \ref{LAX HOLIM}.

\begin{lemma}
\label{pullback etop}
Assume that we have a pullback square
$$\xymatrix{
E' \ar[r]^{\bar{f}}\ar[d]_{p'}& E\ar[d]^{p} \\
B'\ar[r]_{f}& B \\
}$$ 
where $p$, $p'$ are fiber bundles of compact topological manifolds. Let 
$\bar{f}^{\ast}\colon {\Aex}(p) \to {\Aex}(p')$ denote the map induced by the pullback. 
Then there exists a path in ${\Aex}(p)$ joining $\etop(p')$ with $\bar{f}^{\ast}\etop(p)$. 
\end{lemma}

\begin{lemma}
\label{commute}
Assume that we have a commutative diagram 
$$\xymatrix{
E_{1} \ar[r]^{\bar{f}}\ar[d]_{p_{1}}& E_{2}\ar[d]^{p_{2}} \\
B_{1}\ar[r]_{f}& B_{2} \\
E'_{1}\ar[r]_{\bar{f}'}\ar[u]^{p'_{1}} \ar@(l, l)[uu]^{j_{1}}& E'_{2}\ar[u]_{p'_{2}}\ar@(r, r)[uu]_{j_{2}}\\
}$$ 
where for $i=1, 2$ the maps $p_{i}$ are bundles of compact topological manifolds, $p'_{i}$ is 
a sub-bundle of $p_{i}$, $j_{i}\colon E'_{i}\ra E_{i}$ is the inclusion map, and both squares
in the middle are pullbacks.  Then $j_{1\ast}\bar{f}'^{\ast}=\bar{f}^{\ast}j_{2\ast}$, where 
$j_{i\ast}\colon \Aex(p_{i}')\to \Aex(p_{i})$, $\bar{f}^{\ast}\colon \Aex(p_{2})\to \Aex(p_{1})$, 
$\bar{f}'^{\ast}\colon \Aex(p'_{2})\to \Aex(p'_{1})$ are the maps induced by $j_{i}, \bar{f}, \bar{f'}$
respectively.
\end{lemma}

\begin{proof}[Proof of Theorem \ref{MAIN2}]
Let $\fib{}{}$ be a fiber bundle as in the statement of  Theorem \ref{MAIN2}, 
and let $c, c_{i}, \bar{c}, \bar{c}_{i}$ 
denote the maps as in the pullback squares (\ref{UPULL}) above. Applying the map $\bar{c}^{\ast}$ 
to the path $\sigma_{p^{U}}$ from Proposition \ref{universal add} we obtain a path $\bar{c}^{\ast}\sigma_{p^{U}}$ in $\Aex(p)$ joining $\bar{c}^{\ast}\etop(p^{U})$
with $\bar{c}^{\ast}\iota_{1\ast}\etop(p^{U}_{1})+\bar{c}^{\ast}\iota_{2\ast}\etop(p^{U}_{2})-\bar{c}^{\ast}\iota_{0\ast}\etop(p^{U}_{0})$. Lemma \ref{pullback etop}
implies existence of  a path $\eta$ joining $\etop(p)$ with $\bar{c}^{\ast}\etop(p^{U})$. Similarly 
for $i=0, 1, 2$ we have paths $\eta_{i}$ in $\Aex(p_{i})$ joining 
$\bar{c}_{i}^{\ast}\etop(p^{U}_{i})$  with $\etop(p_{i})$. Take 
$\eta':= j_{1\ast}\eta_{1}+ j_{2\ast}\eta_{2}- j_{0\ast}\eta_{0}$. This is a path 
in $\Aex(p)$ with endpoints  
$j_{1\ast}\bar{c}_{1}^{\ast}\etop(p^{U}_{1})+j_{2\ast}\bar{c}_{2}^{\ast}\etop(p^{U}_{2})-
j_{0\ast}\bar{c}_{0}^{\ast}\etop(p^{U}_{0})$ and
$j_{1\ast}\etop(p_{1})+j_{1\ast}\etop(p_{2})-j_{1\ast}\etop(p_{0})$. 
By Lemma \ref{commute} we have $j_{i\ast}\bar{c}_{i}^{\ast}= \bar{c}^{\ast}\iota_{\ast}$. 
Therefore we can concatenate $\eta$, $\bar{c}^{\ast}\sigma_{p^{U}}$ and $\eta'$ and 
obtain a path in $\Aex(p)$ which joins $\etop(p)$ with 
$j_{1\ast} \etop(p_1) +j_{2\ast} \etop(p_2)-j_{0\ast} \etop(p_0)$.
\end{proof}

%%%%%%%%%%%%%%%%%%%%%%%%%%%%%%%%%%%%
%
%                 FROM TOP TO DISCR
%
%%%%%%%%%%%%%%%%%%%%%%%%%%%%%%%%%%%%

\section{From topological to discrete structure group}
\label{FLAT}

As the previous section demonstrated, checking additivity for the parametrized topological Euler 
characteristic reduces to showing that it holds for the universal bundle 
$p^{U}\colon\utot\to B\tops$. The structure group of this bundle is the simplicial group $\tops$. Our 
next goal is to show that it is enough to verify additivity for a certain flat bundle, i.e. a bundle whose structure group is discrete. 

Recall that a map $f\colon B' \ra B$ is a homology equivalence if it induces isomorphisms 
on homology groups with arbitrary local coefficients. 

\begin{lemma}
\label{induced add}
Let $\qfib$ be a fiber bundle with a compact topological manifold $M=M_{1}\cup_{M_{0}}M_{2}$
as a fiber. Assume that $p$ admits a decomposition into sub-bundles as in the statement 
in Theorem \ref{MAIN2}, and that we have a pullback diagram
$$
 \xymatrix{
E'\ar[d]_{p'}\ar[r]^{\bar f}& E\ar[d]^{p} \\
B'\ar[r]^{f} & B \\
}
$$ 
If the map $f\colon B'\ra B$ is a homology equivalence, then the additivity path for $\etop$
exists for the bundle $p$ if it exists for the bundle $p'$ (where $p'$ comes with decomposition 
into sub-bundles induced from $p$). 
\end{lemma}

Lemma \ref{induced add} follows from the following

\begin{proposition}
\label{we sections}
Assume that we have a homotopy pullback square
$$
 \xymatrix{
X'\ar[d]_{q'}\ar[r]^{\bar f}& X\ar[d]^{q} \\
B'\ar[r]^{f} & B \\
}
$$
such that $q$, $q'$ are fibrations whose fibers are nilpotent spaces. Assume also that 
$f\colon B\to B'$ is a homology equivalence. Then the induced map of the spaces of sections
$$\Gamma(q)\to \Gamma(q')$$
is a weak equivalence.  
\end{proposition}

\begin{proof}
See \cite[Proof of Cor. 2.7]{DWW}
\end{proof}

\begin{proof}[Proof of Lemma \ref{induced add}]
Let $p_{i}\colon E_{i}\to B$ be the sub-bundles in the decomposition of $p$, 
and let $j_{i}\colon E_{i}\to E$ be the inclusion maps. Using the assumption 
that we have  an additivity path for $p'$ and Lemma \ref{pullback etop} we obtain a path 
in $A(p')$ joining $\bar f^{\ast}\etop(p)$ with 
$\bar f^{\ast}(j_{1\ast}\etop(p_{1})+j_{2\ast}\etop(p_{2})-j_{0\ast}\etop(p_{0}))$
Thus it suffices to prove that the map $\bar f^{\ast}\colon \Aex(p)\to\Aex(p')$ is a 
weak equivalence. In order to show that consider the homotopy pullback square
$$
 \xymatrix{
\Aex_{B'}(E')\ar[d]_{p'_{\ast}}\ar[r]^{\bar f}& \Aex_{B}(E)\ar[d]^{p_{\ast}} \\
B'\ar[r]^{f} & B \\
}
$$
where the fibrations $p'_{\ast}$, $p_{\ast}$ are defined as in \ref{FAH}.
The fiber of $p_{\ast}$ and $p'_{\ast}$ is $\Aex(M)$ which is an infinite loop space, thus in particular 
a nilpotent space. By Lemma \ref{we sections} we have a weak equivalence of the spaces 
of sections
$$\Gamma(p'_{\ast})\to\Gamma(p_{\ast})$$
On the other hand we have $\Gamma(p_{\ast})\simeq \Aex(p)$, and 
$\Gamma(p'_{\ast})\simeq \Aex(p')$ (see \ref{tb2b}, \ref{FAH}). 
It follows that $\Aex(p)\simeq \Aex(p')$. 
\end{proof}

Our application of Lemma \ref{induced add} is as follows. Consider the universal bundle
$p^{U}\colon \utot\ra B\tops$ as in Section \ref{UNIVERSAL}. By the collaring theorem 
\cite{EK}  the submanifold $M_{0}$ has a bicollar in $M$ i.e. we have an embedding
$c\colon M_{0}\times (-1, 1)\to M$ such 
that $c(m, 0)=m$ for $m\in M_{0}$. Let $T_{c}$ denote the subgroup of $T$ consisting
of all these splitting preserving homeomorphisms of $M$ which are product maps on 
the bicollar $c$.  In other words, $f\in T_{c}$ if there is a homeomorphism 
$f'\colon M_{0}\to M_{0}$ such that the following diagram commutes: 
 $$\xymatrix{
M\ar[r]^{f}& M \\
M_{0}\times(-1, 1)\ar[r]^{f'\times{\rm id}}\ar[u]^{c}& M_{0}\times(-1, 1)\ar[u]_{c} \\
}$$ 
Let $T_{\varepsilon}:={\rm colim_{c}} T_{c}$ where the colimit is taken over all bicollar 
neighborhoods of $M_{0}$, and let  $T^{\delta}_{\varepsilon}$ denote the group $T_{\varepsilon}$, 
but equipped with the discrete topology. We have

\begin{lemma}
\label{TDELTA}
The homomorphism $T^{\delta}_{\varepsilon}\to T$ induces a homology 
equivalence of classifying spaces $BT^{\delta}_{\varepsilon}\to BT$.
\end{lemma}

\begin{proof}
The homomorphism  $T^{\delta}_{\varepsilon}\to T$ is a composition 
of the inclusion $T_{\varepsilon}\to T$ and the map $T^{\delta}_{\varepsilon}\to T_{\varepsilon}$.
The map of classifying spaces induced by the first of these homomorphisms $BT_{\varepsilon}\to BT$
is a homotopy equivalence by Siebenmann's isotopy extension theorem \cite{Siebe}. The map 
$BT^{\delta}_{\varepsilon}\to BT_{\varepsilon}$ induced by the second homomorphism is a 
homology equivalence by the results of McDuff \cite{Mc}, Thurston \cite{Th}, Segal \cite{Sg2}, and Mather \cite{Ma}.
\end{proof}

We are now in position to give a proof of  Lemma \ref{universal add}, and thus complete
the proof of Theorem \ref{MAIN2}.

\begin{proof}[Proof of Lemma \ref{universal add}]
We have a pullback diagram 
$$
 \xymatrix{
ET^{\delta}_{\varepsilon}\times_{T^{\delta}_{\varepsilon}}M \ar[d]_{p^{T^{\delta}_{\varepsilon}}}\ar[r]& \utot\ar[d]^{p^{U}} \\
BT^{\delta}_{\varepsilon}\ar[r] & B\tops \\
}
$$
By Lemma \ref{TDELTA} the map $BT^{\delta}_{\varepsilon}\to BT$ is a homology equivalence. 
Moreover, by Proposition \ref{DISCR BUNDLE ADD} the additivity of $\etop$ holds for the bundle $p^{T^{\delta}_{\varepsilon}}$. Therefore using Lemma \ref{induced add} we obtain
additivity of $\etop$ for the bundle $p^{U}$.
\end{proof}

%%%%%%%%%%%%%%%%%%%%%%%%%%%%%%%%%%%%
%
%                 ADDITIVITY FORMULAS
%
%%%%%%%%%%%%%%%%%%%%%%%%%%%%%%%%%%%%

\section{Additivity for topological Reidemeister torsion}
\label{ADDITIVITY}

Our final task is to give the proof of additivity of the topological Reidemeister torsion, 
i.e. Theorem \ref{MAIN}. 
Let then $\qfib$ be a bundle of compact topological manifolds and let $V\overset{\rho}{\ra}E$ be a  locally constant sheaf of finitely generated projective left $R$-modules such that  the assumptions of Theorem \ref{MAIN} are satisfied. From the pullback diagram below Definition \ref{PARA T TOR}
it follows that in order to prove additivity for $\ttop(p)$ we need to construct 
\begin{itemize}
\item[1)] 
a path $\widetilde\omega^{t}_{p}$ in $\wtA^{\%}(p)$ joining $\wtetop(p)$ with 
$j_{1\ast}\wtetop(p_{1})+j_{2\ast}\wtetop(p_{2})-j_{0\ast}\wtetop(p_{0})$;
 
\item[2)] 
a path $\omega^{h}_{\rho}$ in $\Phi^{h}_{\rho}(p)$ joining $\thom(p)$ with 
$j_{1\ast}\thom(p_{1})+j_{2\ast}\thom(p_{2})-j_{0\ast}\thom(p_{0})$
and such that $\delta_{1}(\widetilde\omega^{\%})=\delta_{2}(\omega^{h}_{\rho})$. 

\end{itemize}  

{\noindent \it Construction of the path $\widetilde\omega^{t}_{p}$.} Recall that the space 
$\wtA^{\%}(p)$ was obtained 
as a homotopy pullback of the diagram  
$$\holim[x]A(E_{x})\overset{\beta}{\lra} \holim[(x, \theta)] A(E_{x}^{\theta})
\overset{\alpha}{\lla}\holim[(x, \theta)]\Aex(E_{x}^{\theta})$$ 
and that the point $\wtetop(p)\in\wtA(p)$ was represented by the triple $(\ehom(p), \sigma_{p}, \etop(p))$
where $\sigma_{p}$ is the canonical path in $\holimin[(x, \theta)] A(E_{x}^{\theta})$ 
joining the images of $\ehom(p)$ and $\etop(p)$. As a consequence in order to describe the path $\widetilde\omega^{t}_{p}$
we need to construct 
\begin{itemize}
\item[(i)] 
a path $\omega^{t}_{p}$ in $\holimin[(x, \theta)] \Aex(E^{\theta}_{x})$ joining $\etop(p)$
with $j_{1\ast}\etop(p_{1})+j_{2\ast}\etop(p_{2})-j_{0\ast}\etop(p_{0})$;
\item[(ii)]
a path $\omega^{h}_{p}$ in $\holimin[x]A(E_{x})$ joining $\ehom(p)$ with 
$j_{1\ast}\ehom(p_{1})+j_{2\ast}\ehom(p_{2})-j_{0\ast}\ehom(p_{0})$;
\item[(iii)] a homotopy $H\colon [0,1]\times[0,1]\to \holimin[(x,\theta)] A(E^{\theta}_{x})$ such that 
$H(-, 0)=\alpha\circ\omega^{t}_{p}$, $H(-, 1)=\beta\circ\omega^{h}_{p}$, $H(0, -)=\sigma_{\rho}$, and
$H(1, -)= j_{1\ast}\sigma_{p_{1}}+j_{2\ast}\sigma_{p_{2}}-j_{0\ast}\sigma_{p_{0}}$. 
\end{itemize}
We take $\omega^{t}_{p}$ to be the additivity path for the topological Euler characteristic which we 
constructed in that proof of Theorem \ref{MAIN2}. In order to construct a suitable path $\omega^{h}_{p}$
we will use $\Abiv(p)$ -- the bivariant $A$-theory of the bundle $p$ \cite{Wi}. The space 
$\Abiv(p)$ is the $K$-theory space of the Waldhausen category $\Rfd(p)$ which is a  full subcategory of the category of retractive spaces over $E$. The objects of $\Rfd(p)$ are these retractive spaces 
$r \colon{X}\rlas{E}\colon s$ which satisfy the condition that $s$ is a cofibrations and that  for every point $b\in B$ the 
homotopy fiber $F^{b}_{p\circ r}$ of the map $p\circ r$ over $b$ is a homotopy finitely dominated space. 
Notice that $F^{b}_{p\circ r}$ is in a natural way a retractive space over $F^{b}_{p}$ -- the homotopy 
fiber of $p$, thus the above assumption says that for every $b\in B$ the fiber $F^{b}_{p\circ r}$
is an object of the Waldhausen category $\Rfd(F^{b}_{p})$ (\ref{NONPARA E}).  
The Waldhausen category structure on $\Rfd(p)$ is defined by taking a morphism to be  
a weak equivalence or a cofibrations if its underlying map of spaces is 
a homotopy equivalence or respectively a cofibration.

 The retractive space ${E\sqcup E}\rlas{E}$ is an object of $\Rfd(p)$, and so it defines a point 
 $\ebiv(p)\in\Abiv(p)$. Let $j_{i\ast}\colon \Abiv(p_{i})\to\Abiv(p)$, $i=0,1,2$, be the maps induced by
 inclusions of subbundles.  Using the same constructions as in the proof of 
 Theorem \ref{nonparam top additivity} with $E$, $E_{i}$ taken in place of $M$, $M_{i}$ we can construct 
 a path $\omega^{biv}_{p}$ in $\Abiv(p)$ joining $\ebiv(p)$ with 
 $j_{\ast1}\ebiv(p_{1})+j_{2\ast}\ebiv(p_{2})-j_{0\ast}\ebiv(p_{0})$.

For every simplex $x\in\simp$ we have an  exact functor $\Rfd(p)\to \Rfd(E_{x})$
which yields a map of infinite loop spaces $\Abiv(p)\to A(E_{x})$. 
These maps can be combined to give a map 
$$\holim[x] \Abiv(p)\to \holim[x] A(E_{x})=A(p)$$
where the first homotopy limit is taken over the constant functor with the value $\Abiv(p)$. Composing this 
map with the Bousfield-Kan map $\Abiv(p)\simeq \lim_{x}\Abiv(p)\to\holimin[x]\Abiv(p)$ we obtain 
$$\alpha^{\ast}\colon \Abiv(p)\to A(p)$$
(this is the generalized coassembly map of \cite{Wi}). The image of $\ebiv(p)$
under the map $\alpha^{\ast}$ does not coincide with $\ehom(p)$, but there is a canonical path 
in $\holimin[x]A(E_{x})$ which joins these two points. As a consequence the path 
$\alpha^{\ast}\circ \omega^{biv}_{p}$ defines an additivity path $\omega^{h}_{p}$ for the homotopy 
Euler characteristic.

In order to see that we have the required homotopy $H$ one needs to retrace our construction 
of the additivity path for the topological Euler characteristic. First, one considers bundles with a 
discrete structure group where the homotopy $H$ can described using constructions on the level 
of Waldhausen categories. For more general bundles one uses the fact that the additivity path 
$\omega^{t}_{p}$ was obtained from the above special case by means of pullbacks. This is enough to 
verify that the homotopy $H$ will still exist. 

{\vskip .2cm \noindent \it Construction of the  path $\omega^{h}_{\rho}$.}
Notice that the path $\omega^{h}_{\rho}$ is an additivity path for the homotopy Reidemeister 
torsion of the bundle $p$. The existence of such path was proved in \cite{Do}. In our present 
setting we need the path $\omega^{h}_{\rho}$ which is compatible with $\widetilde\omega^{t}_{p}$, 
thus it is not enough to quote that result. The path $\omega^{h}_{\rho}$ can be obtained 
however by mimicking the constructions of \cite{Do}.    Briefly, one starts by constructing 
$A^{ac}(p)$ - the acyclic bivariant $A$-theory of the bundle $p$. The space $A^{ac}(p)$
is obtained from a Waldhausen category $\R^{ac}(p)$ which is a full subcategory of 
$\Rfd(p)$. The objects of $\R^{ac}(p)$ are these retractive spaces $\retrlabel{X}{E}$ which satisfy
the condition that for every $b\in B$ the relative chain complex  of the pair $(F_{p\circ r}, F_{r})$
with local coefficients in $\rho|_{F_{r}}$ is acyclic (as before $F_{p\circ r}$ and  $F_{p}$ denote here 
the homotopy fibers over $b$ of the maps $p\circ r$ and $p$ respectively). 
The inclusion of Waldhausen categories $\R^{ac}(p)\hookrightarrow \Rfd(p)$ induces a map
$A^{ac}(p)\to A^{biv}(p)$. Directly from the construction of the path $\omega^{biv}_{p}$ it follows that 
it admits a lift $\omega^{ac}$ to $A^{ac}(p)$. The space  $A^{ac}(p)$ in turn maps into the 
homotopy Whitehead space $\Phi^{h}_{\rho}(p)$ in such way that we obtain a commutative 
diagram 
$$
\xymatrix{
A^{ac}(p)\ar[r]\ar[d] & \Phi^{h}_{\rho}(p)\ar[d] \\
A^{biv}(p)\ar[r] & A(p) \\
}
$$
We can take $\omega^{h}_{\rho}$ to be 
the image of $\omega^{ac}$ under this map.
 
\bibliographystyle{plain}
\bibliography{additivity}

\end{document}